\documentclass[12pt, oneside]{amsart}
\usepackage{amssymb,latexsym}
\usepackage{graphicx}
\usepackage{pdfsync}  
\usepackage[text={6.4in,9in},centering]{geometry}
\usepackage{hyperref}

\newcommand{\arXiv}[2]{\href{https://arxiv.org/abs/#1}{\texttt{arXiv:#1 [#2]}}}

\newcommand{\Q}[2]{Q_{#1}^{(#2)}}
\newcommand{\Asub}[2]{A_{#1}^{(#2)}(t)}
\newcommand{\A}[1]{A^{(#1)}(z)}

\renewcommand{\sf}{.6}

\newcommand{\up}[2]{\raise#1\hbox{#2}}

\usepackage{xcolor}

\begin{document}
\title{A note on Stirling permutations}
\author{Ira M. Gessel}
\address{Department of Mathematics\\
   Brandeis University\\
   Waltham, MA 02453}
\email{gessel@brandeis.edu}

\begin{abstract}
In this note we generalize an identity of John Riordan and Robert Donaghey relating the enumerator for ``Stirling permutations" to the Eulerian polynomials.\end{abstract}

\maketitle
\thispagestyle{empty}

\emph{This paper was written in 1978, but not published until now (May, 2020). I am reproducing the 1978 manuscript here with a few minor corrections.}

\bigskip
For $r\ge1$, let $\Q nr$ be the set of all permutations $a_1a_2\cdots a_{nr}$ of the multiset $\{1^r, 2^r,\dots, n^r\}$ such that if $i<j<k$ and $a_i=a_k$ then $a_j\ge a_i$. For example, a typical permutation in $\Q43$ is
\begin{equation}
\label{e-1}
122211334443.
\end{equation}
It will be convenient to let $\Q n0$ be the one-element set containing the unique increasing permutation of $\{1,2,\dots, n\}$; thus, $\Q 40 = \{1234\}$. We take $\Q 0r$ to be the one-element set containing the ``empty permutation" $\varnothing$. 

More generally, for any set $S$ of positive integers, we let $\Q Sr$ be the analogous set of permutations of $S$, so that $\Q nr = \Q{\{1,2,\dots, n\}}r$. We let $\Q{}r = \bigcup_S \Q Sr$, and we call the elements of $\Q{}r$ \emph{$r$-permutations}. Thus $1$-permutations are ordinary permutations, and $2$-permutations are the ``Stirling permutations" of \cite{3}.

A \emph{descent} (or fall) of a sequence $a_1a_2\cdots a_n$ of integers is an index $i$ for which $a_i>a_{i+1}$. In addition, we count a ``conventional" descent at the end of every \emph{nonempty} sequence. The \emph{descent number} $d(\pi)$ is the number of descents of the sequence $\pi$. (Thus $d(\varnothing) =0$ and $d(\pi)\ge1$ if $\pi\ne\varnothing$.)

Now let
\[\Asub nr = \sum_{\pi\in \Q nr}t^{d(\pi)}\]
and let 
\[\A r = \sum_{n=0}^\infty \Asub nr \frac{z^n}{n!}.\]
(Thus $\Asub n1$ is the ordinary Eulerian polynomial and $\A0 = 1+t(e^z-1)$.)

A differential equation (analogous to the first formula on page 33 of \cite{3}) for $\A r$ is easily obtained: if we remove the 1s from a permutation in $\Q nr$ we obtain a sequence $\pi_1, \pi_2,\dots, \pi_{r+1}$ of $r$-permutations (some of which may be empty). For example, from the 3-permutation \eqref{e-1} we get
\[\varnothing, 222, \varnothing, 334443.\]
The descent number of $\pi$ is the sum of the descent numbers of the $\pi_i$ unless $\pi_{r+1}=\varnothing$, in which case $\pi$ has an additional descent. This decomposition leads to the differential equation
\[\frac{d\ }{dz} \A r = [\A r]^r [\A r -1] +t [\A r]^r.\]

Riordan and Donaghey \cite{4} have found a relationship between $\A 2$ and $\A 1$, which we now generalize to a relationship between $\A r$ and $\A s$ for $0\le s< r$. This relationship is a consequence of a decomposition for $r$-permutations which may be described most easily with the help of a tree representation of $r$-permutations due in the case $r=1$ to Foata and Strehl \cite{1,2}.

We may represent the decomposition described above for the 3-permutation $122211334443$ as
\[ \includegraphics[scale=\sf]{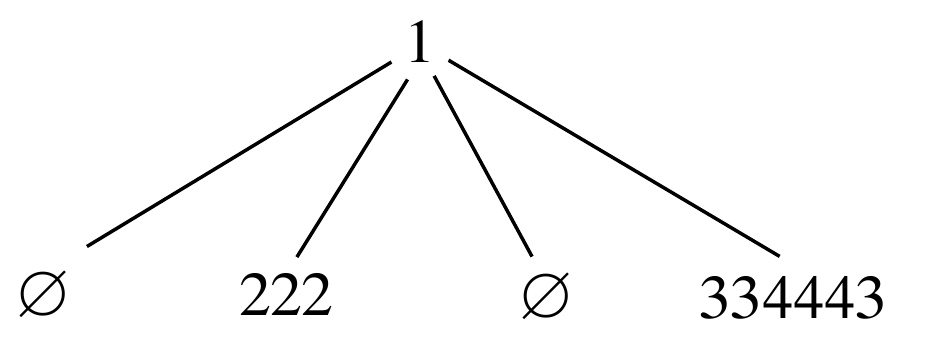}\up{7pt}{\!\!\!.}\]
We now iterate this decomposition to obtain the tree
\begin{equation}
\label{e-2}
\lower 1.3in\hbox{\includegraphics[scale=\sf]{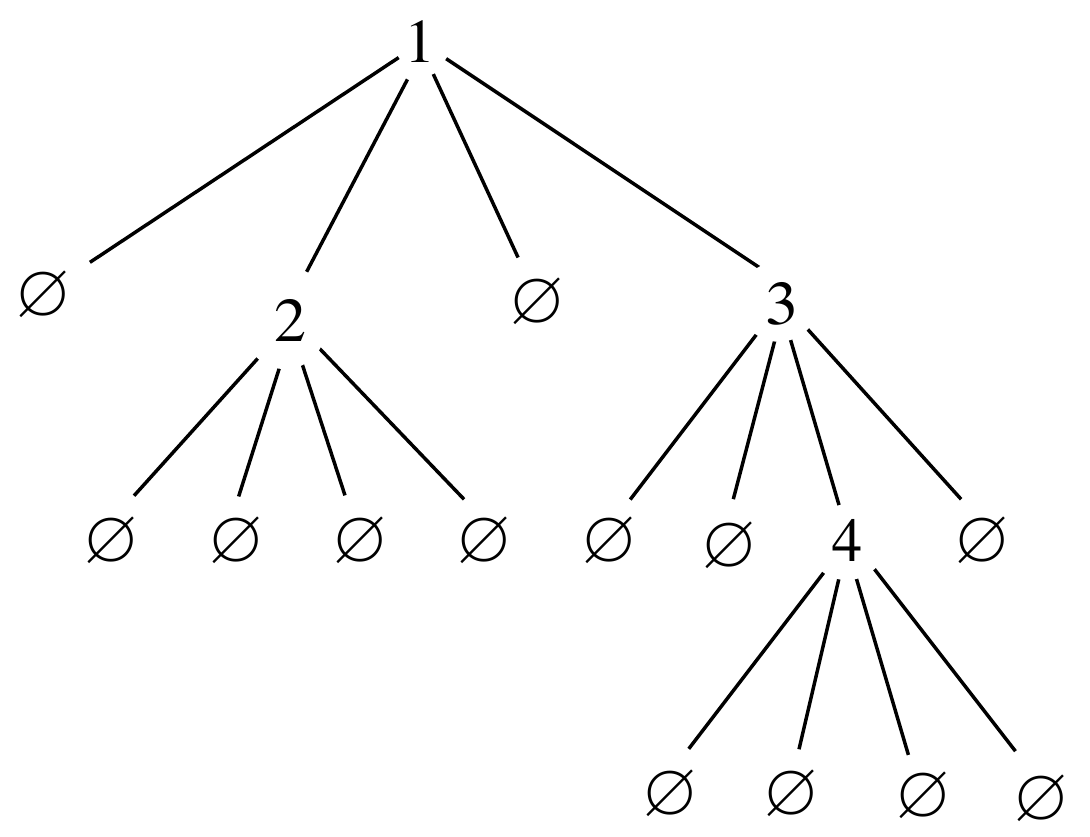}}\up{-91pt}{\!.}
\end{equation}
Thus an $r$-permutation $\pi$ corresponds to an ``$(r+1)$-ary increasing tree" $T(\pi)$, which we call an \emph{$r$-tree}. (A $0$-permutation corresponds to a ``unary tree"---thus $T(123)$ is the tree 
\[\includegraphics[scale=\sf]{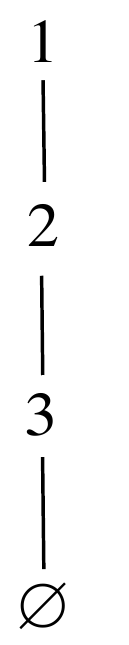}\raise4pt\hbox{\!\!\!\!.)}\]

An important fact is that the descents of $\pi$ are easily read off from $T(\pi)$. Given a node in an $r$-tree, its children are the $r+1$ subtrees lying directly under it, ordered from left to right. Then we observe that only the \emph{last} occurrence of an integer $j$ in an $r$-permutation $\pi$ can be ``followed by" a descent, and that this will happen if and only if the last child of $j$ in $T(\pi)$ is empty. Thus in the 3-permutation $122211334443$ there are descents after $2$, $3$, and $4$, and in the corresponding tree \eqref{e-2}, $2$, $3$, and $4$ have empty last children.

Given an $r$-tree $\tau$ we define its \emph{$s$-skeleton} for $0\le s<r$ to be the $s$-tree obtained from $\tau$ by deleting the first $r-s$ children (and all their descendants) of each node of $\tau$. Thus the 2-skeleton of the tree \eqref{e-2} is
\[\includegraphics[scale=\sf]{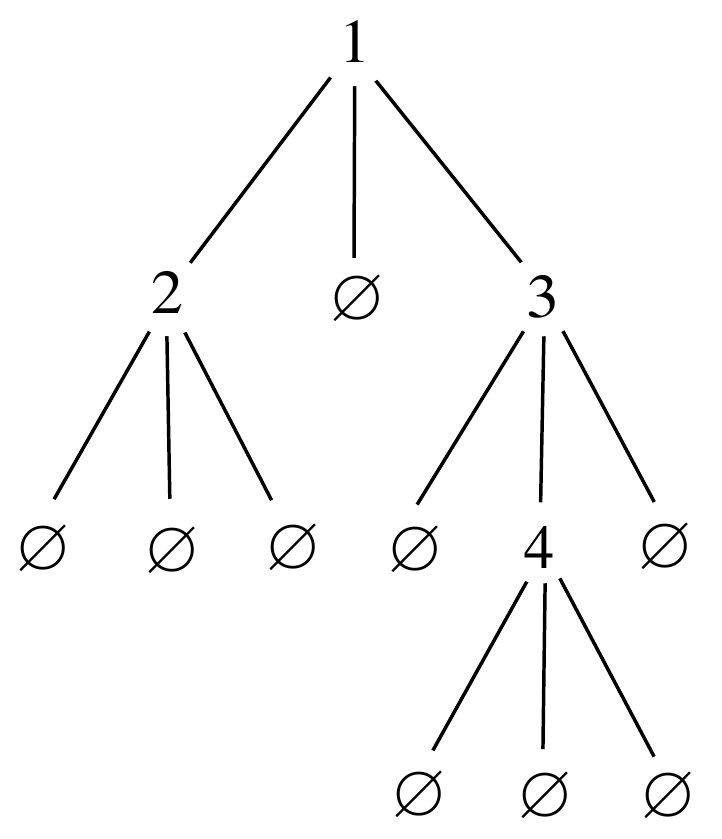}\]
and its $0$-skeleton is 
\[\includegraphics[scale=\sf]{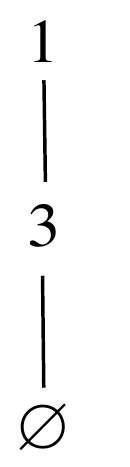}\up{3pt}{\!\!\!\!.}\]
It is easy to see that the number of descents of $\tau$ is equal to the number of descents of the $s$-skeleton and the deleted subtrees. This observation leads to the relationship
\begin{equation*}
\A r = A^{(s)}\left(\int [\A r]^{r-s}\, dz\right).
\end{equation*}
(Here $\int f(z)\, dz$ means $\int_0^z f(u)\, du$, and the outer parentheses on the right indicate functional composition.) In particular, for $s=r-1$ we have
\begin{equation*}
\A r = A^{(r-1)}\left(\int \A r\,dz\right),
\end{equation*}
which for $r=2$ is Riordan and Donaghey's identity.
The tree representation can also be used to generalize many of Foata and Strehl's results \cite{1,2} to $r$-permutations.


%

\bigskip
\centerline{\vrule height .4 pt width 2in}
\bigskip

\textbf{Additional references.} Since this paper was written, a fairly large literature  on Stirling permutations has developed. Here are some of the papers on this topic:

\def\bibitem#1{}
\parindent0pt
\parskip 5pt

\bibitem{MR3414183}
J.~Fernando Barbero~G., Jes\'{u}s Salas, and Eduardo J.~S. Villase\~{n}or,
  \emph{Generalized {S}tirling permutations and forests: higher-order
  {E}ulerian and {W}ard numbers}, Electron. J. Combin. \textbf{22} (2015),
  no.~3, Paper 3.37, 20 pp.

\bibitem{MR2476838}
Mikl\'{o}s B\'{o}na, \emph{Real zeros and normal distribution for statistics on
  {S}tirling permutations defined by {G}essel and {S}tanley}, SIAM J. Discrete
  Math. \textbf{23} (2008/09),  401--406.

\bibitem{MR3555882}
David Callan, Shi-Mei Ma, and Toufik Mansour, \emph{Restricted {S}tirling
  permutations}, Taiwanese J. Math. \textbf{20} (2016), no.~5, 957--978.

\bibitem{MR3828757}
Guan-Huei Duh, Yen-Chi~Roger Lin, Shi-Mei Ma, and Yeong-Nan Yeh, \emph{Some
  statistics on {S}tirling permutations and {S}tirling derangements}, Discrete
  Math. \textbf{341} (2018),  2478--2484.

\bibitem{MR570212}
Dominique Dumont, \emph{Une g\'{e}n\'{e}ralisation trivari\'{e}e sym\'{e}trique
  des nombres eul\'{e}riens}, J. Combin. Theory Ser. A \textbf{28} (1980),
 307--320.

\bibitem{MR3131903}
Askar Dzhumadil'daev and Damir Yeliussizov, \emph{Stirling permutations on
  multisets}, European J. Combin. \textbf{36} (2014), 377--392.

\bibitem{MR3505221}
Rafael~S. Gonz\'{a}lez~D'Le\'{o}n, \emph{On the free {L}ie algebra with
  multiple brackets}, Adv. in Appl. Math. \textbf{79} (2016), 37--97.

\bibitem{MR3983744}
Rafael~S. Gonz\'{a}lez~D'Le\'{o}n, \emph{A family of symmetric functions associated with {S}tirling
  permutations}, J. Comb. \textbf{10} (2019), no.~4, 675--709.

\bibitem{MR2864433}
J.~Haglund and Mirk\'{o} Visontai, \emph{Stable multivariate {E}ulerian
  polynomials and generalized {S}tirling permutations}, European J. Combin.
  \textbf{33} (2012),  477--487.

\bibitem{jkp}
Svante Janson, Markus Kuba, and Alois Panholzer, \emph{Generalized {S}tirling
  permutations, families of increasing trees and urn models}, J. Combin. Theory
  Ser. A \textbf{118} (2011),  94--114.

\bibitem{MR1140465}
Paul Klingsberg and Cynthia Schmalzried, \emph{A family of constructive
  bijections involving {S}tirling permutations}, Proceedings of the
  {T}wenty-first {S}outheastern {C}onference on {C}ombinatorics, {G}raph
  {T}heory, and {C}omputing ({B}oca {R}aton, {FL}, 1990), vol.~78, 1990,
  pp.~11--15.

\bibitem{MR1267287}
Paul Klingsberg and Cynthia Schmalzried,
 \emph{Barred permutations}, Proceedings of the {T}wenty-fourth
  {S}outheastern {I}nternational {C}onference on {C}ombinatorics, {G}raph
  {T}heory, and {C}omputing ({B}oca {R}aton, {FL}, 1993), vol.~95, 1993,
  pp.~153--161.

\bibitem{MR3739896}
Arnold Knopfmacher, Shi-Mei Ma, Toufik Mansour, and Stephan Wagner,
  \emph{Geometrically distributed {S}tirling words and {S}tirling
  compositions}, J. Math. Anal. Appl. \textbf{460} (2018),  98--120.

\bibitem{MR2847273}
Markus Kuba and Alois Panholzer, \emph{Analysis of statistics for generalized
  {S}tirling permutations}, Combin. Probab. Comput. \textbf{20} (2011), 
  875--910.

\bibitem{MR2957938}
Markus Kuba and Alois Panholzer, \emph{Enumeration formul{\ae} for pattern restricted {S}tirling
  permutations}, Discrete Math. \textbf{312} (2012), 3179--3194.

\bibitem{MR3949472}
Markus Kuba and Alois Panholzer, \emph{Stirling permutations containing a single pattern of length
  three}, Australas. J. Combin. \textbf{74} (2019), 215--239.

\bibitem{MR3940782}
Shi-Mei Ma, Jun Ma, and Yeong-Nan Yeh, \emph{The ascent-plateau statistics on
  {S}tirling permutations}, Electron. J. Combin. \textbf{26} (2019), no.~2,
  Paper No. 2.5, 13.

\bibitem{MR4038671}
Shi-Mei Ma, Jun Ma, and Yeong-Nan Yeh, \emph{David-{B}arton type identities and alternating run polynomials},
  Adv. in Appl. Math. \textbf{114} (2020), 101978, 19 pp.

\bibitem{MR1297383}
SeungKyung Park, \emph{Inverse descents of {$r$}-multipermutations}, Discrete
  Math. \textbf{132} (1994),  215--229.

\bibitem{MR1295782}
SeungKyung Park, \emph{{$P$}-partitions and {$q$}-{S}tirling numbers}, J. Combin.
  Theory Ser. A \textbf{68} (1994), 33--52.

\bibitem{MR1280598}
SeungKyung Park, \emph{The {$r$}-multipermutations}, J. Combin. Theory Ser. A
  \textbf{67} (1994), 44--71.

\bibitem{MR3338849}
Jeffrey~B. Remmel and Andrew~Timothy Wilson, \emph{Block patterns in {S}tirling
  permutations}, J. Comb. \textbf{6} (2015), no.~1--2, 179--204.

\bibitem{xu2002generalizations}
Dapeng Xu, \emph{Generalizations of two-stack-sortable permutations}, Ph.D.
  Thesis, Brandeis University, 2002. 
 \arXiv{math/0209313}{math.CO}.


\begin{thebibliography}{9}
\bibitem{1}
  D. Foata and V. Strehl, Rearrangements of the symmetric group and enumerative properties of the tangent and secant numbers, Math. Z. \textbf{137} (1974), 257--264. 
   
\bibitem{2}
D. Foata and V. Strehl, Euler numbers and variations of permutations, Colloquio Internazionale sulle Teorie Combintorie (Roma, 1973), Tom I, 119--131. Atti dei Convegni Lincei, No.~17, Accad. Naz. Lincei, Rome, 1976.

\bibitem{3}
I. Gessel and R. P. Stanley, Stirling polynomials, J. Combinatorial Theory (A) \textbf{24} (1978), 24-33.

\bibitem{4}
J. Riordan, letter to Richard Stanley, April 11, 1978.



\end{thebibliography}

\end{document}